\documentclass[a4paper,reqno,11pt]{amsart}
\usepackage[a4paper,margin=2.5cm]{geometry}
\usepackage[T1]{fontenc}
\usepackage[utf8]{inputenc}
\usepackage[english]{babel}
\usepackage{amssymb,amsmath,amsthm,mathtools}
\usepackage[dvipsnames,svgnames,table]{xcolor}
\usepackage{graphicx}
\usepackage{caption,subcaption}
\usepackage{lmodern}
\usepackage[foot]{amsaddr}
\usepackage[shortlabels]{enumitem}
\usepackage[unicode=true]{hyperref}
\usepackage[capitalise,noabbrev]{cleveref}
\usepackage{thmtools,thm-restate}
\usepackage{comment}
\usepackage{mdframed}
\usepackage{algorithmicx}
\usepackage[noend]{algpseudocode}

\DeclarePairedDelimiter\set{\{}{\}}

\newcommand{\defin}[1]{\emph{\textcolor{ForestGreen}{#1}}}
\newcommand{\defmath}[1]{\defin{\text{$#1$}}}

\newcommand{\Oh}{\mathcal{O}}

\newcommand{\calB}{\mathcal{B}}

\newcommand{\calF}{\mathcal{F}}

\let\le\leqslant
\let\ge\geqslant
\let\leq\leqslant
\let\geq\geqslant

\let\subset\subseteq

\let\epsilon\varepsilon

\let\setminus\backslash

\crefformat{equation}{#2(#1)#3}
\crefformat{subsection}{Subsection #2#1#3}

\makeatletter
\def\thm@space@setup{
  \thm@preskip=4mm
  \thm@postskip=0mm
}
\makeatother

\algdef{SE}[DOWHILE]{Do}{doWhile}{\algorithmicdo}[1]{\algorithmicwhile\ #1}

\mdfdefinestyle{dontsplit}{
  hidealllines=true,
  nobreak=true,
  leftmargin=0pt,
  rightmargin=0pt,
  innerleftmargin=0pt,
  innerrightmargin=0pt,
}

\newmdtheoremenv[style=dontsplit]{theorem}{Theorem}
\newmdtheoremenv[style=dontsplit]{lemma}[theorem]{Lemma}
\newmdtheoremenv[style=dontsplit]{proposition}[theorem]{Proposition}
\newmdtheoremenv[style=dontsplit]{question}[theorem]{Question} 
\newmdtheoremenv[style=dontsplit]{corollary}[theorem]{Corollary} 
\newmdtheoremenv[style=dontsplit]{problem}[theorem]{Problem}
\newmdtheoremenv[style=dontsplit]{conjecture}[theorem]{Conjecture}

\newtheorem*{problem*}{Problem}
\newtheorem*{conjecture*}{Conjecture}
\newtheorem*{remark*}{Remark}
\newtheorem*{question*}{Question} 

\theoremstyle{remark}
\newmdtheoremenv[style=dontsplit]{claim}[theorem]{Claim}

\crefname{claim}{Claim}{Claims}
\newtheorem*{claim*}{Claim} 
\newmdtheoremenv[style=dontsplit]{observation}[theorem]{Observation}

\newenvironment{proofclaim}[1][]
	{\vspace{-\topsep}\begin{proof}[Proof of the claim] }{\end{proof}}

\graphicspath{{figs/}}

\setlist[enumerate,1]{label=\textup{(\roman*)}}

\makeatletter
\newcommand{\myitem}[1]{%
\item[#1]\protected@edef\@currentlabel{#1}%
}
\makeatother

\makeatletter
\expandafter\patchcmd\csname\string\proof\endcsname
  {\topsep6\p@\@plus6\p@\relax}
  {\topsep0pt\relax}
  {}
  {\PackageError{mytemplate}{Could not adjust proof spacing}{}}
\makeatother

\hypersetup{
    colorlinks,
    linkcolor={RoyalBlue},
    citecolor={RubineRed},
    urlcolor={blue!80!black},
    pdftitle={Diagonal Boolean Ramsey numbers are superlinear},
    pdfauthor={Krystian Gajdzica, Jędrzej Hodor, Piotr Micek}
}

\newcommand{\Prob}{\mathbb{P}}
\newcommand{\E}{\mathbb{E}}
\DeclareMathOperator{\Poisson}{Poisson}

\begin{document}

\title[Diagonal Boolean Ramsey numbers are superlinear]
{Diagonal Boolean Ramsey numbers are superlinear}

\author[K.~Gajdzica]{Krystian Gajdzica}
\email{\href{mailto:krystian.gajdzica@uj.edu.pl}{krystian.gajdzica@uj.edu.pl}}

\author[J.~Hodor]{J\k{e}drzej Hodor}
\email{\href{mailto:jedrzej.hodor@gmail.com}{jedrzej.hodor@gmail.com}}

\author[P.~Micek]{Piotr Micek}
\email{\href{mailto:piotr.micek@uj.edu.pl}{piotr.micek@uj.edu.pl}}

\address[K.~Gajdzica, P.~Micek]{Theoretical Computer Science Department, Faculty of Mathematics and Computer Science, Jagiellonian University, Kraków, Poland.}

\address[J.~Hodor]{Theoretical Computer Science Department,
Faculty of Mathematics and Computer Science and Doctoral School of Exact and Natural Sciences,
Jagiellonian University, Krak\'ow, Poland}

\thanks{K. Gajdzica was supported by the National Science Centre, Poland, under grant no.\ 2024/53/N/ST1/01538. The research of J.~Hodor and P.~Micek was partially supported by the National Science Centre, Poland, under grant no.\ 2022/47/B/ST6/02837 within the OPUS 24 program.}

\newcommand{\jedrzej}[1]{\textcolor{ProcessBlue}{Jędrzej: #1}}
\newcommand{\piotr}[1]{\textcolor{orange}{Piotr: #1}}
\newcommand{\krystian}[1]{\textcolor{LimeGreen}{Krystian: #1}}

\begin{abstract}
For a positive integer $n$, the Boolean lattice of order $n$, denoted $\calB_n$, is the family of all subsets of $\{1,\dots,n\}$ ordered by inclusion. 
The diagonal Boolean Ramsey number $R(n)$ 
is the least integer $N$ such that every $2$-coloring of $\calB_N$ contains a monochromatic induced copy of $\calB_n$.
We show that $R(n) = \widetilde{\Omega}(n^{3/2})$, which improves the previous linear bound of Ivan and Wessels. The upper bound remains $\mathcal{O}(n^2)$. 
The proof was found by OpenAI's ChatGPT-6 Astra.
\end{abstract}

\maketitle

\section{Introduction}\label{sec:introduction}

For a positive integer $n$, the \defin{Boolean lattice} of \defin{order} $n$, denoted \defin{$\calB_n$}, is 
the poset $(2^{[n]},\subseteq)$.
For positive integers $N$ and $n$, an \defin{embedding} of $\calB_n$ into $\calB_N$ is an injective map
$\phi\colon 2^{[n]}\to 2^{[N]}$ such that, for all $A,B \in 2^{[n]}$,
$A \subset B$ if and only if $\phi(A)\subset \phi(B)$.
The image of such an embedding is an \defin{induced copy} of $\calB_n$ in $\calB_N$.
The \defin{diagonal Boolean Ramsey number $R(n)$} 
is the least integer $N$ such that every $2$-coloring of $\calB_N$ contains a monochromatic induced copy of $\calB_n$.

It is easy to see that $R(n)\geq 2n$; indeed, consider a coloring of $\calB_{2n-1}$ in which exactly $n$ levels are red and exactly $n$ levels are blue.\footnote{A \defin{level} of a Boolean lattice is a family of all elements of a given cardinality.}
Clearly, there is no monochromatic chain of size $n+1$, so there is no monochromatic $\calB_n$ either.
Also, it is easy to explain that $R(n)=\mathcal{O}(n^2)$. 
The lexicographic product $P[P]$ of any poset $P$ is obtained from $P$ by replacing each element with a copy of $P$.
Kierstead and Trotter~\cite{KT87} showed that every $2$-coloring of $P[P]$ contains a monochromatic copy of $P$. 
Since $\calB_n[\calB_n]$ can be embedded in $\calB_{n^2+2n}$, this proves that $R(n)\leq n^2+2n$. 
This was first observed by Axenovich and Walzer~\cite[Theorem 1.(i)]{AW17}.

It turned out to be surprisingly difficult to improve these bounds for $R(n)$~\cite{CS18,LT22,BP23,Win25}. 
So far, the best known lower bound was given by Ivan and Wessels~\cite{IW26} and the best known upper bound was proved by Axenovich and Winter~\cite[Corollary~2]{AW26}; these bounds are of the following form: 
\[
    2.7n + C < R(n) \leq n^2-(1-o(1))n\log_2 n.
\]
We prove the first superlinear lower bound on $R(n)$.



\begin{theorem}\label{thm:main}
For every sufficiently large integer $n$, 
\[
    R(n)>\frac{1}{2^{13}} \cdot \left(\frac{n}{\ln n}\right)^{3/2}.
\]
\end{theorem}

Although very natural itself, the diagonal Boolean Ramsey numbers are a special case of the following notion introduced by Axenovich and Walzer~\cite{AW17}.
For posets $P$ and $Q$, the \defin{poset Ramsey number} \defin{$R(P,Q)$} is the least integer $N$ such that every red-blue coloring of $\calB_N$ contains a red induced copy of $P$
or a blue induced copy of $Q$. 
It is well-known that every poset can be embedded in a Boolean lattice of large enough order.
Therefore, it is natural to study the Ramsey numbers of $P$ and $Q$ when both of them are Boolean lattices. 


Our lower bound is given by a random $2$-coloring of $\calB_N$ with $N\sim (\frac{n}{\ln n})^{\frac{3}{2}}$, in which, with non-zero probability, there is no monochromatic copy of $\calB_n$. 
We partition the sizes $1,\ldots,N$ into $T=\Theta(\sqrt{n/\ln n})$
intervals of width $W=\Theta(n/\ln n)$, with $N=TW$. Each interval has
$2N$ independent families, each consisting of a Poisson number of random
subsets of $[N]$. The coloring is as follows: a set $S\subseteq[N]$ is blue if it meets every member of
some family assigned to its interval. All other sets are red. 
This test is inspired by a random monotone formula of Talagrand~\cite[proof of
Proposition~5.1]{Tal96}.

We did not try to optimize the polylogarithmic function, nor the multiplicative constant in front of it. The simplicity of the argument appeared to be more important to us, given that there is still a big gap in the exponent between the lower and the upper bound. We leave an obvious open problem.

\begin{problem*}
What is the infimum $\alpha^*$ over real numbers $\alpha$ such that $R(n) = \Oh(n^\alpha)$?
So far we know that $\frac{3}{2}\leq \alpha^*\leq2$.
\end{problem*}

\section{Preliminaries}

For every positive integer $n$, we write $\defmath{[n]} = \{1,\dots,n\}$.
For a finite multiset $\calF$ of subsets of $[N]$, a set $S\subseteq[N]$ is a \defin{transversal} of $\calF$ if $S\cap F\neq\emptyset$ for every $F\in\calF$.


Let $\lambda$ be a positive real number.
We say that a random variable $K$ has a \defin{Poisson distribution} with the parameter $\lambda$ if for every nonnegative integer $k$, $\mathbb{P}(K = k) = e^{-\lambda}\cdot\frac{\lambda^k}{k!}$.
In this case, we write $K \sim \Poisson(\lambda)$.
Our construction can be given without using the Poisson distribution. 
However, certain inequalities are slightly less technical when using Poissons. 
We use standard properties of this distribution. 
First, if $K\sim\Poisson(\lambda)$ for a real number $\lambda$, then $\E[K] = \lambda$.
Second is the following simple fact.
\begin{observation}
\label{obs:Poisson-and-binomial}
    Let $\lambda$ be a positive real number and $p\in(0,1]$.
    Let $K\sim\Poisson(\lambda)$,
    and let $L$ be a random variable.
    Assume that for every nonnegative integer $k$, the distribution of $L$ conditional on $K=k$ is the binomial distribution $\mathrm{Binom}(k,p)$.
    Then $L\sim\Poisson(\lambda p)$.
\end{observation}


\section{The coloring}\label{sec:proof} 

In this section, we prove~\cref{thm:main}. 
To this end, we must construct a certain coloring of a Boolean lattice.
This is done through a random process. 
The random process depends on an input positive integer $n$.
Some events that may or may not occur will be considered.
For such an event, we use a shorthand \defin{with high probability} when the probability of the event tends to one as $n$ approaches infinity.
We prove that the random experiment has the desired properties for the proof of \cref{thm:main} with high probability.
Thus, in many places, we use the fact that $n$ is sufficiently large.
We do not specify explicitly what this means, but it could be inferred from the proofs.

Throughout the paper, we use the following auxiliary parameters that are defined for $n$ large enough

\begin{equation*}\label{eq:parameters}
\begin{gathered}
    r=\left\lceil\sqrt{n\ln n}\right\rceil,
    \qquad
    W=\left\lfloor\frac{n}{2048\ln n}\right\rfloor,
    \qquad
    T=\left\lfloor\frac{n}{2r}\right\rfloor,
    \\
    N=TW,
    \qquad
    p=\frac{1}{64W},
    \qquad
    d=\left\lfloor\frac{np}{4}\right\rfloor.
\end{gathered}
\end{equation*}

Let $(K_{j,i} : j\in [T], i \in [2N])$ be a family of independent random variables such that for all $j \in [T]$ and $i \in [2N]$, we have
\begin{equation*}\label{eq:lambda}
    K_{j,i}\sim\Poisson(\lambda_j) \qquad \text{and}
    \qquad
    \lambda_j=\frac r4(1-p)^{-jW}.
\end{equation*}

Let $\calF_{j,i}$ be a multiset of $K_{j,i}$ subsets of $[N]$ chosen independently in the following way.
We include each element of $[N]$ in the set with probability $p$ independently of all other choices.

Given a set $S\subseteq [N]$ and $j\in[T]$, we define
\[
Z_j(S)=|\set{i\in[2N]: \text{$S$ is a transversal of $\calF_{j,i}$}}|.
\]
The following is a straightforward observation:
\begin{equation} \label{eq:mono:ZjS}
    \text{ for all $S,S'\subset [N]$ and $j \in [T]$, if $S \subset S'$ and $Z_j(S) > 0$, then $Z_j(S') > 0$.}
\end{equation}

Depending on the described random experiment, we define a \textcolor{red}{red}-\textcolor{blue}{blue} coloring of $\calB_N$. 
The empty set is \textcolor{red}{red}.
Let $S\subseteq [N]$ be nonempty, and let $j \in [T]$ be the unique integer such that $(j-1)W<|S|\le jW$.
The set $S$ is
\begin{align*}
    \text{\textcolor{blue}{blue} }& \text{if }  Z_j(S)>0 \text{ and }\\
    \text{\textcolor{red}{red} }& \text{if } Z_j(S)=0.
\end{align*}

The \defin{output of the random experiment} (or just the \defin{output} for short) is a consistent assignment of the values of all considered random variables, and in particular, the \textcolor{red}{red}-\textcolor{blue}{blue} coloring of $\calB_N$.
We establish three properties of the output that hold with high probability.
First, we state them (see \ref{prop:add}--\ref{prop:chain}) and show that they are sufficient to complete the proof of \cref{thm:main}, more precisely, to show that with high probability $\calB_N$ does not contain a monochromatic copy of $\calB_n$.
Then, in \cref{sec:props}, we prove the properties.

\begin{enumerate}[label=(p\arabic*)]
    \item\label{prop:add}
    For every $j\in[T]$ and all disjoint sets $P,Y\subseteq[N]$ with
    $|P|\geq(j-1)W$ and $|Y|\geq n/2$, there exists a set $D\subseteq Y$ of size $r$
    such that $Z_j(P\cup D)>0$.
    \item\label{prop:delete}
    For every $j\in[T]$ and all sets $A\subseteq Q\subseteq[N]$ with
    $|Q|\leq jW$ and $|A|\geq n/2$, there exists a set $D\subseteq A$ such that
    $|D|\leq d\cdot Z_j(Q)$ and $Z_j(Q\setminus D)=0$.
    \item\label{prop:chain}
    For every sequence $Q_1\subseteq\cdots\subseteq Q_T\subseteq[N]$ with
    $|Q_j|\leq jW$ for every $j\in[T]$,
    \[
        d\sum_{j=1}^T Z_j(Q_j)<n/2.
    \]
\end{enumerate}

\begin{proof}[Proof of \cref{thm:main} given an output of the experiment satisfying \ref{prop:add}--\ref{prop:chain}]
We prove that if an output satisfies \ref{prop:add}--\ref{prop:chain}, then $\calB_N$ does not contain a monochromatic copy of $\calB_n$.
This is indeed sufficient to obtain \cref{thm:main}.
Fix the values of all random variables from the experiment as in the given output.

Before we begin the main part of the proof, we state a simple observation that was first noted by Axenovich and Walzer~\cite[Theorem~9]{AW17}.
Given an embedding $\phi$ of $\calB_n$ into $\calB_N$, for every $x \in [n]$, one may choose a \emph{private witness} $w_x \in \phi(\{x\}) \setminus \phi([n] \setminus\{x\})$.
Then, with $X = \{w_x : x \in [n]\}$, we have that for every $A \subset X$, $\phi(A) \cap X = \{w_x : x \in A\}$.
For simplicity, given an embedding of $\calB_n$ into $\calB_N$, we may rename the elements of $[N]$ (and all dependent objects) so that $w_x = x$ for every $x \in [n]$; then also $X = [n]$.
For future reference, we restate the property of $\phi$ that may or may not hold:
\begin{equation}\label{eq:trace}
    \text{for every $A \subset [n]$}, \qquad \phi(A) \cap [n] = A.
\end{equation}

First suppose that there is a \textcolor{red}{red} induced copy of $\calB_n$ in $\calB_N$ embedded by a map $\phi$.
As discussed before, we may assume that \cref{eq:trace} holds.
We construct sets $A_1\subseteq \cdots\subseteq A_T\subseteq [n]$ such that
\[ |A_j|=jr \quad \text{and} \quad |\phi(A_j)|>jW \qquad \text{ for every $j\in[T]$.}\]
Note that this leads to a contradiction as $\phi(A_T)$ is supposed to have more elements than $TW = N$.
The construction is by induction.
We add to the construction $A_0=\emptyset$.

Let $j\in[T]$ and suppose that $A_{j-1}$ has been constructed.
Set $P=\phi(A_{j-1})$ and $Y=[n]\setminus A_{j-1}$.
By \cref{eq:trace}, $P\cap Y=\emptyset$.
By the induction hypothesis, $|P|\geq(j-1)W$ and by the definition of $r$ and $T$, we have $|Y|=n-(j-1)r\geq n-rT\geq n/2$.
Therefore, we may apply \ref{prop:add}, obtaining $D\subseteq Y$ of size $r$ such that $Z_j(P\cup D) > 0$.

Let $A_j=A_{j-1}\cup D$.
Since $\phi$ is inclusion-preserving, we have $P\cup D\subseteq\phi(A_j)$.
In particular, by \cref{eq:mono:ZjS}, $Z_j(\phi(A_j)) > 0$.
Also, $P\cap D=\emptyset$ and $r\geq1$, so $|\phi(A_j)|\geq|P|+r>(j-1)W$.
Since the copy is \textcolor{red}{red}, the set $\phi(A_j)$ is red.
This forces $|\phi(A_j)|>jW$, which is exactly what is required for the induction to continue.
As discussed, this is a contradiction showing that there is no \textcolor{red}{red} induced copy of $\calB_n$ in $\calB_N$.

Next suppose that there is a \textcolor{blue}{blue} induced copy of $\calB_n$ in $\calB_N$ embedded by a map $\phi$.
Again, we may assume that \cref{eq:trace} holds.
We construct sets $A_0\subseteq \cdots\subseteq A_T\subseteq [n]$ such that
\[ |A_j| > n/2 \quad \text{and} \quad |\phi(A_j)|\leq jW \qquad \text{ for every $j\in[T] \cup \{0\}$.}\] 
Note that this leads to a contradiction as $\phi(A_0)$ is supposed to have zero elements, which must be false as $|A_0| > 0$.
The construction is by reverse induction.
We begin with $A_T = [n]$.

Let $j\in[T]$ and suppose that the sets $A_{j},\dots,A_T$ satisfying the required conditions have been constructed.
Note that $A_j\subseteq \phi(A_j)$ by \cref{eq:trace}.
Property~\ref{prop:delete} yields a set $D_j\subseteq A_j$ such that $|D_j|\leq d\cdot Z_j(\phi(A_j))$ and $Z_j(\phi(A_j)\setminus D_j)=0$.

Let $A_{j-1}=A_j\setminus D_j$.
Observe that $\phi(A_{j-1}) \subset \phi(A_j) \setminus D_j$ and by \cref{eq:mono:ZjS}, $Z_j(\phi(A_{j-1})) = 0$.
Since the copy is \textcolor{blue}{blue}, the set $\phi(A_{j-1})$ is blue and $|\phi(A_{j-1})| \leq |\phi(A_j)| \leq jW$.
This forces $|\phi(A_{j-1})| \leq (j-1)W$, which is one of the inductive conditions that we must preserve.

It remains to verify that $|A_{j-1}| > n/2$.
The sets $\phi(A_j)\subseteq\cdots\subseteq \phi(A_T)$ form a nested sequence with
$|\phi(A_i)|\leq iW$ for every $i\in\{j,\dots,T\}$.
Extend this sequence to indices $1,\dots,T$ by placing $\emptyset$ at every earlier index.
Property~\ref{prop:chain} applies to this sequence giving
\[
    d\sum_{i=j}^T Z_i(\phi(A_i))<n/2.
\]
Consequently,
\begin{equation*}\label{eq:blue-budget}
    |A_{j-1}|
       =n-\sum_{i=j}^T|D_i|
       \geq n-d\sum_{i=j}^T Z_i(\phi(A_i))>n/2.
\end{equation*}
Thus, the second induction invariant follows.
This completes the inductive proof.
As discussed, this is a contradiction showing that there is no \textcolor{blue}{blue} induced copy of $\calB_n$ in $\calB_N$.
This ends the proof of the theorem.
\end{proof}

\section{The properties of the coloring}\label{sec:props}


\begin{claim}
    \label{claim:add}
    With high probability \ref{prop:add} holds.
\end{claim}

\begin{proofclaim}
Fix $j\in[T]$ and disjoint sets $P,Y\subseteq[N]$ satisfying the assumptions of \ref{prop:add}.
For a fixed $i\in[2N]$, let $U_i$ be a random variable counting the members of $\calF_{j,i}$ that are disjoint from $P$,
and let $V_i$ be a random variable counting the members of $\calF_{j,i}$ that are disjoint from $P\cup Y$.
By~\cref{obs:Poisson-and-binomial},  
$U_i\sim\Poisson(\lambda_j(1-p)^{|P|})$ and 
$V_i\sim\Poisson(\lambda_j(1-p)^{|P|+|Y|})$. 
Thus, for sufficiently large $n$,
\[
    \E U_i=\lambda_j(1-p)^{|P|}
       =\frac r4(1-p)^{|P|-jW}
       \leq\frac r4(1-p)^{-W} < \frac r3,
\]
where the last inequality follows since $p=1/(64W)$ and $W\to\infty$ imply $(1-p)^{-W}\to e^{1/64}<4/3$.

Since the sets $P$ and $Y$ are disjoint, for sufficiently large $n$, we also have
\begin{align*}
    \E V_i=\lambda_j(1-p)^{|P|+|Y|}
         &=(\E U_i)\cdot (1-p)^{|Y|}\\
         &<r/3 \cdot e^{-pn/2}&&\textrm{by $1-p\leq e^{-p}$ and $|Y|\geq n/2$}\\  
         &\leq r/3 \cdot  n^{-16}&&\textrm{since $pn/2=n/(128W)\ge16\ln n$}\\
         &=o(1)&&\textrm{since $r\sim \sqrt{n\ln n}$.}
\end{align*}
Upper bounds on $\E U_i$ and $\E V_i$ allow us to apply Markov's inequality to obtain the following bound for sufficiently large values of $n$:
\begin{align*}
    \Prob(U_i>r\text{ or }V_i>0) 
    &\leq \Prob(U_i>r) + \Prob(V_i>0)\\
    &\leq\frac{\E U_i}{r}+\E V_i\\
    &\leq\frac13+o(1)<\frac12.
\end{align*}
Thus, with probability at least $1/2$, at most $r$ members of $\calF_{j,i}$ are disjoint from $P$,
and every such member intersects $Y$.
Choose one element of $Y$ from each of these members.
The set of chosen elements has size at most $r$; extend it to an $r$-element set $D\subseteq Y$,
which is possible because $|Y|\geq n/2\geq r$.
Every member of $\calF_{j,i}$ then intersects $P\cup D$.

The random variables $(\calF_{j,i} : i\in[2N])$ are independent.
Hence, the probability that none supplies such a set $D$ is at most $2^{-2N}$. 
There are at most $3^N$ ordered pairs of disjoint sets $P,Y\subseteq[N]$ as 
each element of $[N]$ belongs to $P$, to $Y$, or to neither.
A union bound over these pairs and $j\in[T]$ bounds the failure probability of the event in the statement of the claim by $3^N \cdot 2^{-2N}\cdot T=(3/4)^N\cdot T$.
By the definitions of $N$ and $T$ this value tends to $0$ when $n$ approaches infinity.
\end{proofclaim}

\begin{claim}
    \label{claim:delete}
    With high probability \ref{prop:delete} holds.
\end{claim}

\begin{proofclaim}
Fix $j \in [T]$ and sets $A,Q \subset [N]$ satisfying the assumptions of \ref{prop:delete}.

Let $i \in [2N]$.
We say that $F \subset [N]$ is \defin{$(A,Q,d)$-good} if 
\begin{equation*}\label{eq:witness}
    F\cap Q\subseteq A
    \quad\text{and}\quad
    |F\cap Q|=d.
\end{equation*}
Note that for a member $F$ of $\calF_{j,i}$, we have
\[
   \mathbb{P}(F \text{ is $(A,Q,d)$-good}) = \binom {|A|}{d} p^d(1-p)^{|Q|-d}.
\]
Indeed, we choose the $d$ elements in $A$ that belong to $F$, require that the other elements
of $Q$ do not belong to $F$, and impose no condition outside $Q$.
\cref{obs:Poisson-and-binomial} guarantees that the number of $(A,Q,d)$-good members of $\calF_{j,i}$ is a Poisson random variable with the parameter 
\begin{align*}
    \mu
       &=\lambda_j\binom{|A|}{d} p^d(1-p)^{|Q|-d}=\frac r4(1-p)^{|Q|-jW}\binom{|A|}{d}
           \left(\frac p{1-p}\right)^d\\
       &\geq\frac r4\binom{|A|}{d}
           \left(\frac p{1-p}\right)^d
        \geq \frac r4\left(\frac{|A|}{d}\right)^d p^d \\
        &\geq\frac r4\left(\frac{np}{2d}\right)^d &&\textrm{as $|A| \geq n/2$}\\
        &\geq\frac r4\,2^d &&\textrm{as $d \leq np/4$}.
\end{align*}

Since $d \sim 8 \ln n$ and $r\sim \sqrt{n\ln{n}}$, it follows that
\[
    \frac{\ln((r/4)2^d)}{\ln n}
       =\frac{\ln{r}}{\ln n}-\frac{\ln{4}}{\ln{n}}+\frac{d\ln2}{\ln n}
       \longrightarrow\frac12+8\ln2>4.
\]
In particular, $\mu\geq n^4$ for sufficiently large $n$.
By the definition of the Poisson distribution, the probability that $\calF_{j,i}$ has no $(A,Q,d)$-good members is at most $e^{-\mu} \leq e^{-n^4}$.

Now, there are at most $3^N$ pairs $A\subseteq Q\subseteq[N]$, as each element lies in
$A$, in $Q\setminus A$, or outside $Q$.
The union bound yields that the probability that there exist $(j,i) \in [T] \times [2N]$ and $A \subset Q \subset [N]$ satisfying the assumptions of the claim such that $\calF_{j,i}$ has no $(A,Q,d)$-good member is at most $T\cdot 2N \cdot 3^N\cdot e^{-n^4}$. Since $N=\Theta((n/\ln{n})^{3/2})$ and $T=\Theta((n/\ln{n})^{1/2})$, that value tends to $0$ when $n$ approaches infinity.

Therefore, we assume that for all $(j,i) \in [T] \times [2N]$ and $A \subset Q \subset [N]$ satisfying the assumption of the claim the multiset $\calF_{j,i}$ has a $(A,Q,d)$-good member.
Let 
\[I=\{i\in[2N]:\text{ $Q$ \text{ is a transversal of } $\calF_{j,i}$ }\}.\]
For each $i\in I$, choose one $(A,Q,d)$-good member $F_i\in\calF_{j,i}$, and set
\[\textstyle
    D=\bigcup_{i\in I}(F_i\cap Q).
\]
Then $D\subseteq A$ and $|D|\leq d\cdot |I|=d\cdot Z_j(Q)$.
For every $i\in I$, the set $Q\setminus D$ is disjoint from $F_i$.
For every $i\notin I$, some member of $\calF_{j,i}$ is already disjoint from $Q$ and
remains disjoint from $Q\setminus D$.
Therefore, we get $Z_j(Q\setminus D)=0$, as required.
\end{proofclaim}

\begin{claim}
    \label{claim:chain}
    With high probability \ref{prop:chain} holds.
\end{claim}
\begin{proofclaim}
Let $j\in[T]$, $i\in[2N]$, and $Q\subseteq[N]$ with $|Q|\leq jW$.
A random member of $\calF_{j,i}$ is disjoint from $Q$ with probability $(1-p)^{|Q|}$.
The set $Q$ is a transversal of $\calF_{j,i}$ exactly when no member of $\calF_{j,i}$ is disjoint from it.
Therefore, \cref{obs:Poisson-and-binomial} and the definition of the Poisson distribution ensure that
\begin{align*}\label{eq:transversal-prob}
    \Prob(Q\text{ is a transversal of }\calF_{j,i})
       &=\exp\bigl(-\lambda_j(1-p)^{|Q|}\bigr)\\
       &=\exp\left(-\frac r4(1-p)^{|Q|-jW}\right)
        \leq e^{-r/4}.
\end{align*}
Set $k=\lceil n/(2d)\rceil$.
For a fixed nested sequence $Q_1\subseteq\cdots\subseteq Q_T$ with $|Q_j|\leq jW$ for every $j \in [T]$,
the inequality $\sum_{j=1}^T Z_j(Q_j)\geq k$ implies that there are $k$ distinct pairs
$(j,i)\in[T]\times[2N]$ such that $Q_j$ is a transversal of $\calF_{j,i}$.
For a fixed choice of these pairs, the corresponding events are independent. 
Their joint probability is at most $e^{-rk/4}$.

There are at most $(T+1)^N$ nested sequences of subsets of $[N]$ of length $T$.
To see this, record for each element the first set $Q_j$ containing it, or record that
it appears in none of the sets; this gives $T+1$ possibilities per element.
The union bound yields that the probability that there exists 
a sequence $Q_1\subseteq\cdots\subseteq Q_T\subseteq[N]$ with
    $|Q_j|\leq jW$ for every $j\in[T]$ such that $\sum_{j=1}^T Z_j(Q_j)\geq k$ is at most 
\begin{equation*}
    (T+1)^N\binom{T \cdot 2N}{k}e^{-rk/4}.
\end{equation*}
For sufficiently large $n$, we have $k\leq 2NT$ and
$\ln(2NT)=O(\ln n)=o(r)$, so $2NT\leq e^{r/8}$.
Using $\binom{2NT}{k}\leq (2NT)^k$, we obtain that
\begin{equation}\label{eq:chain-exponent}
    (T+1)^N\binom{T \cdot 2N}{k}e^{-rk/4}\leq\exp\left(N\ln(T+1)-\frac{rk}{8}\right).
\end{equation}
Now, it is enough to verify that the expression on the right hand side of the above inequality tends to zero with $n$ approaching infinity. To do this, let us note that
\begin{equation*}
    N\ln(T+1)\sim\frac{n^{3/2}}{8192\sqrt{\ln n}},
    \qquad
    \frac{rk}{8}\sim\frac{n^{3/2}}{128\sqrt{\ln n}}.
\end{equation*}
Indeed, $\ln(T+1)\sim(\ln n)/2$ and $k\sim n/(16\ln n)$.
It follows that $N\ln(T+1)-rk/8\to-\infty$ as $n\to\infty$, and so $\exp\left(N\ln(T+1)-rk/8\right)$ tends to zero for $n$ approaching infinity, as required.
\end{proofclaim}

Note that if each of the three properties \ref{prop:add}--\ref{prop:chain} holds with high probability, then all of them hold simultaneously with high probability.
Therefore, for $n$ large enough, there exists an output of the random experiment satisfying \ref{prop:add}--\ref{prop:chain}.
This completes the proof of \cref{thm:main}.

\section*{AI statement}
The proof was found by OpenAI’s ChatGPT-6 Astra. 
The authors digested the initially
generated output and made an effort to present the argument in a readable and optimized form.
The authors take responsibility for the mathematical correctness of the presented arguments.

\bibliographystyle{plain}
\bibliography{bibliography}

\end{document}